# Differentiable equivalence of fractional linear maps

**Fritz Schweiger**[1]

*University of Salzburg*

**Abstract:** A Moebius system is an ergodic fibred system $(B, T)$ (see [5]) defined on an interval $B = [a, b]$ with partition $(J_k), k \in I, \#I \geq 2$ such that $Tx = \frac{c_k + d_k x}{a_k + b_k x}$, $x \in J_k$ and $T|_{J_k}$ is a bijective map from $J_k$ onto $B$. It is well known that for $\#I = 2$ the invariant density can be written in the form $h(x) = \int_{B^*} \frac{dy}{(1+xy)^2}$ where $B^*$ is a suitable interval. This result does not hold for $\#I \geq 3$. However, in this paper for $\#I = 3$ two classes of interval maps are determined which allow the extension of the before mentioned result.

## 1. Introduction

**Definition 1.** Let $B$ be an interval and $T : B \to B$ be a map. We assume that there is a countable collection of intervals $(J_k), k \in I, \#I \geq 2$ and an associated sequence of matrices

$$\alpha(k) = \begin{pmatrix} a_k & b_k \\ c_k & d_k \end{pmatrix}$$

where $\det \alpha(k) = a_k d_k - b_k c_k \neq 0$, with the properties:

- $\bigcup_{k \in I} \overline{J_k} = \overline{B}$, $J_m \cap J_n = \emptyset$ if $n \neq m$.
- $Tx = \frac{c_k + d_k x}{a_k + b_k x}$, $x \in J_k$
- $T|_{J_k}$ is a bijective map from $J_k$ onto $B$.

Then we call $(B, T)$ a *Moebius system*.

Examples of Moebius systems are abundant. We mention the $g$-adic map $Tx = gx$ mod 1, $g \geq 2$, $g \in \mathbf{N}$ or the map related with regular continued fractions $Tx = \frac{1}{x}$ mod 1. Another important example is the Rényi map

$$T : [0, 1] \to [0, 1]$$

$$Tx = \frac{x}{1-x}, 0 \leq x \leq \frac{1}{2}; \quad Tx = \frac{1-x}{x}, \frac{1}{2} \leq x \leq 1.$$

A Moebius system is a special case of a fibred system [5]. Since $T|_{J_k}$ is bijective the inverse map $V_k : B \to J_k$ exists. The corresponding matrix will be denoted by

$$\beta(k) = \begin{pmatrix} d_k & -b_k \\ -c_k & a_k \end{pmatrix}.$$

---

[1]Dept. of Mathematics, University of Salzburg, Hellbrunnerstr. 34, A-5020 Salzburg, Austria, e-mail: fritz.schweiger@sbg.ac.at

We denote furthermore $\omega(k;x) := |V_k'(x)| = \frac{|a_k d_k - b_k c_k|}{(d_k - b_k x)^2}$. Then a nonnegative measurable function $h$ is the density of an invariant measure if and only if the Kuzmin equation

$$h(x) = \sum_{k \in I} h(V_k x) \omega(k;x)$$

is satisfied.

**Remark.** It is easy to see that we can assume $B = [a,b]$, $B = [a,\infty[$ or $B = ]-\infty, b]$ but $B = \mathbf{R}$ is excluded: Since $\#I \geq 2$ there must exist a Moebius map from a subinterval onto $B = \mathbf{R}$ which is impossible.

**Definition 2.** A Moebius system $(B^*, T^*)$ is called a *natural dual* of $(B, T)$ if there is a partition $\{J_k^*\}$, $k \in I$, of $B^*$ such that $T^* y = \frac{b_k + d_k y}{a_k + c_k y}$, $y \in J_k^*$ i.e. the matrix $\alpha^*(k)$ is the transposed matrix of $\alpha(k)$.

**Theorem 1.** *If $(B^*, T^*)$ is a natural dual of $(B, T)$ then the density of the invariant measure for the map $T$ is given as $h(x) = \int_{B^*} \frac{dy}{(1+xy)^2}$.*

*Proof.* Starting with the Kuzmin equation this follows from the relation

$$(1 + \frac{-c_k + a_k x}{d_k - b_k x} y)^2 (d_k - b_k x)^2 = (1 + \frac{-b_k + a_k y}{d_k - c_k y} x)^2 (d_k - c_k y)^2.$$

Details are given in [5].    □

**Definition 3.** The Moebius system $(B^*, T^*)$ is *differentiably isomorphic* to $(B, T)$ if there is a bijective map $\psi : B \to B^*$ such that $\psi'$ exists almost everywhere and the commutativity condition $\psi \circ T = T^* \circ \psi$ holds. As an example we consider the Rényi map. The system

$$T^* : [0, \infty[ \to [0, \infty[$$

$$T^* y = \frac{1-y}{y}, 0 < y \leq 1; T^* y = -1 + y, 1 \leq y$$

is a natural dual which is differentiably isomorphic under $\psi(t) = \frac{1-t}{t}$ to the Rényi map.

**Theorem 2.** *If the natural dual system $(B^*, T^*)$ is differentiably isomorphic to $(B, T)$ then $\psi(t) = \frac{b+dt}{a+bt}$.*

*Proof.* Assume for simplicity $B = [0,1]$. Then

$$h(x) = \int_{B^*} \frac{dy}{(1+xy)^2} = \int_0^1 \frac{|\psi'(z)|}{(1+x\psi(z))^2} dz = |\frac{\psi(1)}{1+x\psi(1)} - \frac{\psi(0)}{1+x\psi(0)}|.$$

Since the invariant density for $T^*$ is given as $h^*(y) = \frac{1}{1+y}$ we also find $h(x) = \frac{|\psi'(x)|}{1+\psi(x)}$. Integration gives

$$\psi(x) = \frac{\psi(0) + x(\psi(1) + \psi(1)\psi(0) - \psi(0))}{1 + x\psi(0)}.$$

   □

**Remark.** Note that if the natural dual system $(B^*, T^*)$ is differentiably isomorphic to $(B, T)$ with a map $\psi(t) = \frac{c+dt}{a+bt}$ then it follows from the proof of Theorem 2 that $b = c$ is satisfied.

**Lemma 1.** *Such a map $\psi : B \to B^*$, $\psi(t) = \frac{b+dt}{a+bt}$ exists if and only if the conditions*

$$ab_k + b(d_k - a_k) - dc_k = 0, k \in I$$

*are satisfied.*

*Proof.* From $\psi \circ T = T^* \circ \psi$ we get the equations

$$\begin{pmatrix} a & b \\ b & d \end{pmatrix} \begin{pmatrix} a_k & b_k \\ c_k & d_k \end{pmatrix} = \rho \begin{pmatrix} a_k & c_k \\ b_k & d_k \end{pmatrix} \begin{pmatrix} a & b \\ b & d \end{pmatrix}$$

with a constant $\rho \neq 0$. If $aa_k + bc_k \neq 0$ or $bb_k + dd_k \neq 0$ then $\rho = 1$ and the equation $ab_k + bd_k = a_k b + c_k d$ remains. If $aa_k + bc_k = 0$ and $bb_k + dd_k = 0$ then we see that $\rho^2 = 1$. Only the case $\rho = -1$ needs to be considered. We obtain the equations

$$aa_k + bc_k = 0$$
$$bb_k + dd_k = 0$$
$$ab_k + b(a_k + d_k) + dc_k = 0.$$

A non-trivial solution exists if

$$(a_k + d_k)(b_k c_k - a_k d_k) = 0.$$

Since $\det \alpha(k) \neq 0$ we obtain $a_k + d_k = 0$. But then $\alpha(k)^2 = (a_k^2 + b_k c_k)\mathbf{1}$. This means that $T^2 x = x$ on $J_k$ which is not possible. Therefore only $\rho = 1$ remains. $\square$

**Remark.** Let $I = \{1, 2\}$ then the system of two equations

$$ab_k + b(d_k - a_k) - dc_k = 0, k = 1, 2$$

is always soluble. This explains the result given in [4]. Note that the degenerate case that $B^*$ reduces to a single point carrying Dirac measure is included. This happens for $Tx = 2x, 0 \leq x < \frac{1}{2}; Tx = 2x - 1, \frac{1}{2} \leq x < 1$. Then formally we obtain two branches $T^*y = 2y$ and $T^*y = \frac{2y}{1-y}$. Here $B^* = \{0\}$ and $h(x) = 1$.

**Remark.** Haas [1] constructs invariant measures for a family of Moebius systems. It is easy to verify that in all cases $(B, T)$ has a natural dual $(B^*, T^*)$ which explains the invariant densities given under corollary 2 in [1].

**Definition 4.** A natural dual system $(B^*, T^*)$ which is not differentiably isomorphic to $(B, T)$ is called an *exceptional dual*. We will first give an example of an exceptional dual. We consider a special case of Nakada's continued fractions [3]. Setting

$$g = \frac{-1 + \sqrt{5}}{2}, k = [|\frac{1}{x}| + 1 - g]$$

$$B = [g - 1, g]; Tx = \frac{1}{x} - k, x > 0; Tx = -\frac{1}{x} - k, x < 0$$

$$B^* = [0, \frac{1}{2}]; T^*y = \frac{1}{y} - k, \frac{2}{2k+1} < y < \frac{1}{k}; T^*y = -\frac{1}{y} + k, \frac{1}{k} < y < \frac{2}{2k-1}$$

then we consider the equation

$$\psi \circ T_2 \circ T_{-3} \circ T_{-4} = T_2^* \circ T_{-3}^* \circ T_{-4}^* \circ \psi.$$

However, the equation

$$\begin{pmatrix} a & b \\ b & d \end{pmatrix} \begin{pmatrix} 18 & 5 \\ 11 & 3 \end{pmatrix} = \rho \begin{pmatrix} 26 & -7 \\ 11 & -3 \end{pmatrix} \begin{pmatrix} a & b \\ b & d \end{pmatrix}, |\rho| = 1$$

leads to $a = b = d = 0$.

In this paper we will consider in more detail the case $I = \{1, 2, 3\}$. To avoid complications with parameters we fix the partition as $0 < \frac{1}{2} < \frac{2}{3} < 1$. Note that a different partition would give a Moebius system which is not isomorphic by a Moebius map $\psi(t) = \frac{c+dt}{a+bt}$.

## 2. Differentiably isomorphic dual systems

It is easy to see that on every subinterval $[0, \frac{1}{2}], [\frac{1}{2}, \frac{2}{3}], [\frac{2}{3}, 1]$ we can choose a fractional linear map depending on one parameter $\lambda, \mu, \nu$ say. Furthermore, $T$ restricted to one of these subintervals can be increasing or decreasing. We call $T$ of type $(\epsilon_1, \epsilon_2, \epsilon_3)$ where $\epsilon_j = 1$ stands for "increasing" and $\epsilon_j = -1$ for "decreasing". The parameters satisfy the equations

$$\epsilon_1 \det \alpha(1) = \lambda, \epsilon_2 \det \alpha(2) = \mu, \epsilon_3 \det \alpha(3) = \nu.$$

By the choice of the parameters $\lambda, \mu, \nu$ we have $\lambda > 0$, $\mu > 0$, and $\nu > 0$ but due to the fact that no attractive fixed point is allowed some additional restrictions hold e. g. $0 < \lambda \leq 1$ if $\epsilon_1 = 1$ or $1 \leq \nu$ if $\epsilon_3 = 1$.

Now let us assume that a natural dual $(B^*, T^*)$ exists. Then the branches of $T^*$ corresponding to the parameters $\lambda, \mu, \nu$ could appear in six possible orders from left to right. More precisely, if $J_\lambda, J_\mu, J_\nu$ are the intervals such that $T^*$ is defined piecewise by the matrices $\alpha_1^*, \alpha_2^*, \alpha_3^*$ then the three intervals $J_\lambda, J_\mu, J_\nu$ could be arranged in six possible orders, namely $\lambda\mu\nu, \nu\mu\lambda, \lambda\nu\mu, \mu\nu\lambda, \mu\lambda\nu, \nu\lambda\mu$.

We give a list of the matrices $\alpha$, $\beta$, and $\beta^*$.

$\epsilon_1 = 1$
$$\begin{pmatrix} \lambda & 1-2\lambda \\ 0 & 1 \end{pmatrix} \quad \begin{pmatrix} 1 & 2\lambda-1 \\ 0 & \lambda \end{pmatrix} \quad \begin{pmatrix} 1 & 0 \\ 2\lambda-1 & \lambda \end{pmatrix}$$
$\epsilon_1 = -1$
$$\begin{pmatrix} -1 & -\lambda+2 \\ -1 & 2 \end{pmatrix} \quad \begin{pmatrix} 2 & \lambda-2 \\ 1 & -1 \end{pmatrix} \quad \begin{pmatrix} 2 & 1 \\ \lambda-2 & -1 \end{pmatrix}$$
$\epsilon_2 = 1$
$$\begin{pmatrix} 2\mu-1 & 2-3\mu \\ -1 & 2 \end{pmatrix} \quad \begin{pmatrix} 2 & 3\mu-2 \\ 1 & 2\mu-1 \end{pmatrix} \quad \begin{pmatrix} 2 & 1 \\ 3\mu-2 & 2\mu-1 \end{pmatrix}$$
$\epsilon_2 = -1$
$$\begin{pmatrix} \mu-2 & 3-2\mu \\ -2 & 3 \end{pmatrix} \quad \begin{pmatrix} 3 & 2\mu-3 \\ 2 & \mu-2 \end{pmatrix} \quad \begin{pmatrix} 3 & 2 \\ 2\mu-3 & \mu-2 \end{pmatrix}$$
$\epsilon_3 = 1$
$$\begin{pmatrix} \nu-2 & -\nu+3 \\ -2 & 3 \end{pmatrix} \quad \begin{pmatrix} 3 & \nu-3 \\ 2 & \nu-2 \end{pmatrix} \quad \begin{pmatrix} 3 & 2 \\ \nu-3 & \nu-2 \end{pmatrix}$$
$\epsilon_3 = -1$
$$\begin{pmatrix} 2\nu-1 & 1-3\nu \\ -1 & 1 \end{pmatrix} \quad \begin{pmatrix} 1 & 3\nu-1 \\ 1 & 2\nu-1 \end{pmatrix} \quad \begin{pmatrix} 1 & 1 \\ 3\nu-1 & 2\nu-1 \end{pmatrix}$$

**Lemma 2.** *The natural dual $(B^*, T^*)$ is differentiably isomorphic to $(B, T)$ if and only if the following condition $\mathcal{C}$ holds.*

$$\begin{vmatrix} b_1 & d_1 - a_1 & c_1 \\ b_2 & d_2 - a_2 & c_2 \\ b_3 & d_3 - a_3 & c_3 \end{vmatrix} = 0$$

*Proof.* This is clear from Lemma 1. □

**Lemma 3.** *If $(B^*, T^*)$ is differentiably isomorphic to $(B, T)$, then only the orders $\lambda\mu\nu, \nu\mu\lambda$ can appear.*

*Proof.* Since $\psi(t) = \frac{b+dt}{a+bt}$ is either order preserving or order reversing the assertion is immediate. □

**Theorem 3.** *If the natural dual $(B^*, T^*)$ has the order $\lambda\mu\nu$ or $\nu\mu\lambda$, then $(B^*, T^*)$ is differentiably isomorphic to $(B, T)$. In other words: No exceptional dual exists with orders $\lambda\mu\nu$ or $\nu\mu\lambda$.*

*Proof.* The proof considers all types $(\epsilon_1, \epsilon_2, \epsilon_3)$. It first lists the form of condition $\mathcal{C}$ and then the 'boundary conditions' which means that the three maps of $T^*$ fit together to form a Moebius system.

$\boxed{\text{Type } (1,1,1)}$

$$2\lambda\mu + 2\mu = \lambda\nu + \lambda$$

$$V_\lambda^* \xi = \xi, V_\mu^* \xi = V_\lambda^* \sigma, V_\nu^* \xi = V_\mu^* \sigma, \sigma = V_\nu^* \sigma.$$

Then we find $2\lambda - 1 + \lambda\xi = \xi$ which gives $2\lambda - 1 = \xi(1-\lambda)$. Furthermore

$$\sigma = \frac{\nu - 3 + (\nu-2)\sigma}{3 + 2\sigma}$$

$$2\sigma^2 + \sigma(-\nu + 5) - \nu + 3 = 0.$$

If $\sigma = -1$ then $V_\mu^* \sigma = \mu - 1$, $V_\lambda^* \sigma = \lambda - 1$, $V_\mu^* \xi = \lambda\mu + \mu - 1$, $V_\nu^* \xi = \frac{\lambda\nu - 1 - \lambda}{\lambda + 1}$. Hence we obtain the equations $\lambda\mu + \mu = \lambda$, $\lambda\mu + \mu = \lambda\nu$ which show $\nu = 1$ and $2\lambda\mu + 2\mu = \lambda\nu + \lambda$.
Therefore we concentrate on $\sigma = \frac{\nu - 3}{2}$.
Then we calculate $V_\mu^* \xi = \mu + \lambda\mu - 1, V_\lambda^* \sigma = \frac{\lambda + \lambda\nu - 2}{2}$ from which we again arrive at $2\lambda\mu + 2\mu = \lambda\nu + \lambda$.

$\boxed{\text{Type } (1,1,-1)}$

$$\lambda\mu + \mu = \lambda\nu + \lambda + \nu$$

$$V_\lambda^* \xi = \xi, V_\mu^* \xi = V_\lambda^* \beta, V_\nu^* \beta = V_\mu^* \beta, V_\nu^* \xi = \beta.$$

We find again $2\lambda - 1 = \xi(1-\lambda)$. Then we calculate $V_\mu^* \xi = \mu + \lambda\mu - 1$, $V_\nu^* \xi = \frac{\lambda\nu + \nu - \lambda}{\lambda} = \beta, V_\lambda^* \beta = \lambda\nu + \lambda + \nu - 1$ and we see that $\mu\lambda + \mu - 1 = \nu\lambda + \lambda + \nu - 1$ which is condition $\mathcal{C}$.

$\boxed{\text{Type } (1,-1,-1)}$

$$2\lambda\nu + \lambda + 2\nu = \mu$$

$$V_\lambda^* \xi = \xi, V_\mu^* \xi = V_\nu^* \beta, V_\nu^* \xi = \beta, V_\lambda^* \beta = V_\mu^* \beta.$$

We find again $2\lambda - 1 = \xi(1-\lambda)$. Further we calculate $V_\nu^* \xi = \frac{\lambda\nu + \nu - \lambda}{\lambda} = \beta, V_\nu^* \beta = \frac{2\lambda\nu + 2\nu - 1}{\lambda + 1}, V_\mu^* \xi = \frac{\mu - \lambda - 1}{\lambda + 1}$. This shows $\mu - \lambda = 2\lambda\nu + 2\nu$.

**Type $(1, -1, 1)$**

$$\lambda\nu = \mu$$

$$V_\lambda^*\xi = \xi, V_\mu^*\xi = V_\nu^*\xi, V_\lambda^*\sigma = V_\mu^*\sigma, \sigma = V_\nu^*\sigma.$$

As before $2\lambda - 1 = \xi(1-\lambda)$ and further calculations show $V_\mu^*\xi = \frac{-\lambda+\mu-1}{\lambda+1}$, $V_\nu^*\xi = \frac{-\lambda+\lambda\nu-1}{\lambda+1}$ hence $\lambda\nu = \mu$. It is easy to see that $V_\lambda^*\sigma = V_\mu^*\sigma$ gives the same condition. Note that $\sigma = -1$ corresponds to $\nu = 1$.

**Type $(-1, 1, -1)$**

$$2\lambda\mu + \mu = 2\lambda\nu + \lambda$$

$$V_\lambda^*\beta = \alpha, V_\mu^*\alpha = V_\lambda^*\alpha, V_\mu^*\beta = V_\nu^*\beta, \beta = V_\nu^*\alpha.$$

Then $\alpha = V_\lambda^* V_\nu^* \alpha$ and we obtain $\alpha^2(1+2\nu) + \alpha(5\nu + 2 - \lambda) + 3\nu - \lambda + 1 = 0$. Therefore $\alpha = -1$ or $\alpha = \frac{-3\nu+\lambda-1}{1+2\nu}$. If $\alpha = \frac{-3\nu+\lambda-1}{1+2\nu}$ then $V_\lambda^*\alpha = \frac{-1+2\lambda\nu-\nu}{1+\lambda+\nu}$ and $V_\mu^*\alpha = \frac{-1-\nu-\lambda+\mu+2\mu\lambda}{1+\lambda+\nu}$ which gives immediately condition $\mathcal{C}$. If $\alpha = -1$ then $\beta = \infty$ and an easy calculation gives $V_\lambda^*\alpha = \lambda-1$, $V_\mu^*\alpha = \mu-1$, $V_\nu^*\beta = 2\nu-1$, $V_\mu^*\beta = 2\mu-1$ which shows $\lambda = \mu = \nu$.

**Type $(-1, 1, 1)$**

$$4\mu\lambda + \mu\nu + \mu = \lambda\nu + \nu$$

$$V_\lambda^*\sigma = \alpha, V_\mu^*\alpha = V_\lambda^*\alpha, V_\mu^*\sigma = V_\nu^*\alpha, V_\nu^*\sigma = \sigma.$$

As before $\sigma = -1$ or $\sigma = \frac{\nu-3}{2}$. If $\sigma = -1$ then $\alpha = \lambda - 1$, $V_\lambda^*\alpha = \frac{-1}{1+\lambda}$, $V_\mu^*\alpha = \frac{\mu-1-\lambda+2\mu\lambda}{1+\lambda}$, $V_\nu^*\alpha = \frac{\lambda\nu-2\lambda-1}{1+2\lambda}$, $V_\mu^*\sigma = \mu - 1$. Hence $\mu + 2\lambda\mu = \lambda$, $\mu + 2\lambda\mu = \lambda\nu$ which implies $\nu = 1$ and $\mathcal{C}$ is satisfied.
If $\sigma = \frac{\nu-3}{2}$ then calculation gives

$$\alpha = \frac{2\lambda - \nu - 1}{\nu + 1}, V_\mu^*\alpha = \frac{4\lambda\mu + \mu\nu + \mu}{2\lambda + \nu + 1} - 1, V_\lambda^*\alpha = \frac{\lambda\nu + \lambda}{2\lambda + \nu + 1} - 1,$$

$$V_\nu^*\alpha = \frac{2\lambda\nu}{4\lambda + \nu + 1} - 1, V_\mu^*\sigma = \frac{2\mu\nu}{\nu + 1} - 1$$

hence $\lambda\nu + \lambda = \mu\nu + 4\lambda\mu + \mu$.

**Type $(-1, -1, 1)$**

$$2\lambda\nu = 2\lambda\mu + \mu\nu + \mu$$

$$V_\lambda^*\sigma = \alpha, V_\mu^*\sigma = V_\lambda^*\alpha, V_\mu^*\alpha = V_\nu^*\alpha, V_\nu^*\sigma = \sigma.$$

If $\sigma = -1$ then $\alpha = \lambda - 1$. Calculation shows $V_\lambda^*\alpha = \frac{-1}{1+\lambda}$, $V_\mu^*\sigma = \mu - 1$, $V_\mu^*\alpha = \frac{\lambda\mu-2\lambda+\mu-1}{2\lambda+1}$, $V_\nu^*\alpha = \frac{\lambda\nu-2\lambda-1}{2\lambda+1}$. This gives $\lambda = \mu + \lambda\mu$, $\lambda\nu = \mu + \lambda\mu$ and $\nu = 1$ as expected. Condition $\mathcal{C}$ is satisfied.
Now suppose $\sigma = \frac{\nu-3}{2}$ the a similar calculation shows again that condition $\mathcal{C}$ is satisfied.

Type $(-1,-1,-1)$

$$4\lambda\nu + \lambda + \nu = \mu\nu + \lambda\mu + \mu$$

$$V_\lambda^*\beta = \alpha, V_\mu^*\beta = V_\lambda^*\alpha, V_\mu^*\alpha = V_\nu^*\beta, V_\nu^*\alpha = \beta.$$

Since $\alpha = V_\lambda^*V_\nu^*\alpha$, we obtain

$$\alpha^2(1+2\nu) + \alpha(5\nu - \lambda + 2) + 3\nu - \lambda + 1 = 0.$$

If $\alpha = -1$ then $\beta = \infty, V_\lambda^*\alpha = \lambda - 1, V_\mu^*\beta = \frac{\mu-2}{2}$, $V_\mu^*\alpha = \mu - 1, V_\nu^*\beta = 2\nu - 1$. therefore $\nu = \lambda$ and $\mu = 2\lambda$. Condition $\mathcal{C}$ is satisfied.
Now let $\alpha = \frac{\lambda - 3\nu - 1}{2\nu + 1}$. then $(\lambda - \nu)\beta = 2\nu + 2\nu\lambda - \lambda$. Furthermore we obtain

$$V_\lambda^*\alpha = -1 + \frac{2\lambda\nu + \lambda}{1 + \lambda + \nu}, V_\mu^*\beta = -1 + \frac{\lambda\mu + 2\lambda\mu\nu}{\lambda + \nu + 4\lambda\nu}$$

$$V_\mu^*\alpha = \frac{-1 + \mu + \mu\nu + \lambda\mu - 2\lambda}{2\lambda + 1}, V_\nu^*\beta = \frac{\nu + 4\lambda\nu - \lambda - 1}{2\lambda + 1}$$

From $V_\mu^*\beta = V_\lambda^*\alpha$ we obtain $\frac{1}{1+\lambda+\nu} = \frac{\mu}{\lambda+\nu+4\lambda\nu}$ and from $V_\mu^*\alpha = V_\nu^*\beta$ we get $-1 + \mu + \mu\nu + \lambda\mu - 2\lambda = \nu + 4\lambda\nu - \lambda - 1$. Both equations lead to condition $\mathcal{C}$. $\square$

## 2. Exceptional dual systems

Exceptional dual systems exist for some orders. A full discussion of all possible cases is not only lengthy but also of limited value since the method is similar in all cases. Therefore we give some examples for existence and non-existence of such systems.

Type $(1,1,-1)$

Orders $\lambda\nu\mu$ and $\mu\nu\lambda$

$$\xi = V_\lambda^*\xi, V_\lambda^*\gamma = V_\nu^*\gamma, V_\nu^*\xi = V_\mu^*\xi, \gamma = V_\mu^*\gamma$$

We find $\xi(1 - \lambda) = 2\lambda - 1$ as before. Then $V_\nu^*\xi = \frac{\nu + \lambda\nu - \lambda}{\lambda}$, $V_\mu^*\xi = \mu\lambda + \mu - 1$. This gives the equation $\nu + \nu\lambda = \mu\lambda^2 + \mu\lambda$ and eventually $\nu = \mu\lambda$. From $V_\mu^*\gamma = \gamma$ we find that

$$\gamma^2 + \gamma(3 - 2\mu) - 3\mu + 2 = 0.$$

We use $V_\lambda^*\gamma = 2\lambda - 1 + \lambda\gamma$ and $V_\lambda^* = \frac{3\nu - 1 + (2\nu - 1)\gamma}{1 + \gamma}$ and find the equation

$$\gamma^2\lambda + \gamma(3\lambda - 2\nu) + 2\lambda - 3\nu = 0.$$

Comparing both equations we find again $\nu = \lambda\mu$.

Type $(1,1,-1)$

Orders $\mu\lambda\nu$ and $\nu\lambda\mu$

$$\gamma = V_\mu^*\gamma, V_\lambda^*\gamma = V_\mu^*\beta, V_\lambda^*\beta = V_\nu^*\beta, V_\nu^*\gamma = \beta$$

Again we obtain
$$\gamma^2 + \gamma(3 - 2\mu) - 3\mu + 2 = 0.$$
Therefore $\gamma \neq -1$. From $V_\nu^* \gamma = \beta$ we obtain $\beta = -1 + \frac{3\nu + 2\nu\gamma}{1+\gamma}$.
The equation $V_\lambda^* \gamma = V_\mu^* \beta$ gives $\beta(2\mu - 2\lambda - \lambda\gamma) = -3\mu + 4\lambda + 2\lambda\gamma$. If $2\mu - 2\lambda - \lambda\gamma = 0$ then $-3\mu + 4\lambda + 2\lambda\gamma = 0$ which gives $\mu = 0$, a contradiction.
Therefore
$$(3\nu - 1 + \gamma(2\nu - 1))(2\mu - 2\lambda - \lambda\gamma) = (1 + \gamma)(-3\mu + 4\lambda + 2\lambda\gamma)$$
which gives the equation
$$\gamma^2(\lambda + 2\lambda\nu) + \gamma(3\lambda - \mu + 7\lambda\nu - 4\mu\nu) - \mu + 2\lambda - 6\mu\nu + 6\lambda\nu = 0.$$
Therefore
$$3 - 2\mu = \frac{3\lambda - \mu + 7\lambda\nu - 4\mu\nu}{\lambda + 2\lambda\nu}, -3\mu + 2 = \frac{-\mu + 2\lambda - 6\mu\nu + 6\lambda\nu}{\lambda + 2\lambda\nu}.$$
From this we get
$$4\lambda\mu\nu + \lambda\nu - 4\mu\nu + 2\lambda\mu - \mu = 0, 6\lambda\mu\nu + 2\lambda\nu - 6\mu\nu + 3\lambda\mu - \mu = 0.$$
Eliminating $\lambda\mu\nu$ we obtain $\lambda\nu + \mu = 0$ which is a contradiction since $\lambda, \mu, \nu > 0$.
Therefore no exceptional system exists with this configuration.
If one uses $V_\lambda^* \beta = V_\nu^* \beta$ one gets the equation
$$\beta^2 \lambda + \beta(3\lambda - 2\nu) + 2\lambda - 3\nu = 0.$$
If we insert $\gamma = -1 + \frac{\nu}{\beta - 2\nu + 1}$ in the quadratic equation for $\gamma$ we obtain a second equation for $\beta$ which leads again to $\mu + \nu\lambda = 0$.

Let us give another example of an exceptional dual system with an unforeseen complicated condition.

$\boxed{\text{Type } (1, 1, 1)}$

Orders $\lambda\nu\mu$ and $\mu\nu\lambda$
$$\xi = V_\lambda^* \xi, V_\nu^* \xi = V_\lambda^* \gamma, V_\nu^* \gamma = V_\mu^* \xi, \gamma = V_\mu^* \gamma$$
We find again $\xi(1 - \lambda) = 2\lambda - 1$ and $\gamma^2 + \gamma(3 - 2\mu) - 3\mu + 2 = 0$. Note that $\gamma \neq \frac{-3}{2}$.
Calculation gives
$$V_\nu^* \xi = \frac{\lambda\nu - \lambda - 1}{\lambda + 1}, V_\lambda^* \gamma = 2\lambda - 1 + \lambda\gamma$$
$$V_\mu^* \xi = \lambda\mu + \mu - 1, V_\nu^* \gamma = \frac{\nu - 3 + \gamma(\nu - 2)}{3 + 2\gamma}.$$
Equating $V_\nu^* \xi = V_\lambda^* \gamma$ gives $\gamma = -1 + \frac{\nu - \lambda - 1}{\lambda + 1}$ and $V_\nu^* \gamma = V_\mu^* \xi$ gives $\gamma = -1 - \frac{\lambda\mu + \mu}{2\mu\lambda + 2\mu - \nu}$. From this the condition
$$2\lambda\mu\nu + 2\mu\nu + \lambda\nu + \nu = \nu^2 + 2\mu\lambda + \mu\lambda^2 + \mu$$
follows. If we use $(\gamma + 1)^2 + (\gamma + 1)(1 - 2\mu) - \mu = 0$ and substitute $\gamma + 1 = \frac{\nu - \lambda - 1}{\lambda + 1}$ we get the same condition. As the example $\lambda = \frac{1}{2}, \mu = \frac{4}{15}, \nu = 2$ shows this condition can be satisfied easily.

### Type $(1, 1, 1)$

Orders $\nu\lambda\mu$ and $\mu\lambda\nu$

$$\sigma = V_\nu^* \sigma, V_\lambda^* \sigma = V_\nu^* \gamma, V_\mu^* \sigma = V_\lambda^* \gamma, \gamma = V_\mu^* \gamma$$

We get $\gamma^2 + \gamma(3 - 2\mu) - 3\mu + 2 = 0$. From $V_\lambda^* \sigma = V_\nu^* \gamma$ and $V_\mu^* \sigma = V_\lambda^* \gamma$ we calculate

$$\sigma = -2 + \frac{\nu + \nu\gamma}{3\lambda + 2\lambda\gamma}, \sigma = -2 + \frac{\mu}{2\mu - 2\lambda - \lambda\gamma}$$

and eventually we get

$$\gamma^2 \lambda\nu + \gamma(3\lambda\nu + 2\lambda\mu - 2\mu\nu) + 3\lambda\mu - 2\mu\nu + 2\lambda\nu = 0.$$

Therefore comparing the coefficients of the quadratic equations for $\gamma$ we get $3 - 2\mu = 3 + \frac{2\lambda\mu - 2\mu\nu}{\lambda\nu}$, $2 - 3\mu = 2 + \frac{3\lambda\mu - 2\mu\nu}{\lambda\nu}$ which shows $\mu\nu = 0$, a contradiction.

### Type $(1, -1, -1)$

Orders $\nu\lambda\mu$ and $\mu\lambda\nu$

$$\delta = V_\mu^* \beta, V_\lambda^* \delta = V_\mu^* \delta, V_\nu^* \beta = V_\lambda^* \beta, \beta = V_\nu^* \delta$$

Therefore $V_\nu^* V_\mu^* \beta = \beta$ which gives the quadratic equation

$$\beta^2 \mu + \beta(3\mu - 2\nu\mu - 2\nu) - 4\nu\mu - 3\nu + 2\mu = 0.$$

The equation $V_\nu^* \beta = V_\lambda^* \beta$ gives

$$\beta^2 \lambda + \beta(3\lambda - 2\nu) + 2\lambda - 3\nu = 0.$$

Therefore

$$\frac{3\mu - 2\nu\mu - 2\nu}{\mu} = \frac{3\lambda - 2\nu}{\lambda}, \frac{-4\nu\mu - 3\nu + 2\mu}{\mu} = \frac{2\lambda - 3\nu}{\lambda}$$

and then

$$\nu\mu\lambda + \nu\lambda = \mu\nu, 4\nu\mu\lambda + 3\nu\lambda = 3\mu\nu.$$

Hence $\nu\mu\lambda = 0$, a contradiction.

### Type $(1, -1, -1)$

Orders $\lambda\nu\mu$ and $\mu\nu\lambda$

$$\xi = V_\lambda^* \xi, V_\lambda^* \delta = V_\nu^* \delta, V_\nu^* \xi = V_\mu^* \delta, \delta = V_\mu^* \xi$$

We note $\xi(1 - \lambda) = 2\lambda - 1$, $V_\nu^* \xi = -1 + \frac{\nu + \lambda\nu}{\lambda}$, $V_\mu^* \xi = \delta = -1 + \frac{\mu}{\lambda + 1}$. The equation $V_\lambda^* \delta = V_\nu^* \delta$ gives

$$\delta^2 \lambda + \delta(3\lambda - 2\nu) + 2\lambda - 3\nu = 0.$$

If one inserts $\delta = -1 + \frac{\mu}{\lambda + 1}$ we get

$$\mu^2 \lambda + \mu\lambda^2 + \mu\lambda = 2\nu\mu\lambda + \nu\lambda^2 + 2\nu\mu + 2\lambda\nu + \nu.$$

We further calculate $V_\mu^* \delta = -1 + \frac{\lambda\mu + \mu^2 + \mu}{\lambda + 2\mu + 1}$. The equation $V_\nu^* \xi = V_\mu^* \delta = -1 + \frac{\nu + \nu\lambda}{\lambda}$ gives the same condition. Therefore exceptional systems exist. An example is given by $\lambda = 1, \mu = 4, \nu = \frac{6}{5}$.

$\boxed{\text{Type } (1, -1, 1)}$

Orders $\lambda\nu\mu$ and $\mu\nu\lambda$

$$\xi = V_\lambda^* \xi, V_\lambda^* \delta = V_\nu^* \xi, V_\nu^* \delta = V_\mu^* \delta, \delta = V_\mu^* \xi$$

Calculation gives $\xi(1 - \lambda) = 2\lambda - 1$, $\delta = -1 + \frac{\mu}{\lambda+1}$. The condition $V_\nu^* \delta = V_\mu^* \delta$ gives $\nu = \lambda + \mu + 1$ but the condition $V_\lambda^* \delta = V_\nu^* \xi$ also leads to $\lambda^2 + \lambda\mu + \lambda = \lambda\nu$. Since $\lambda > 0$ this is an equivalent condition. Examples are easy to find ($\lambda = 1, \mu = 1, \nu = 3$).

$\boxed{\text{Type } (1, -1, 1)}$

Orders $\mu\lambda\nu$ and $\nu\lambda\mu$

$$V_\mu^* \delta = V_\lambda^* \delta, V_\mu^* \sigma = \delta, V_\nu^* \delta = V_\lambda^* \sigma, \sigma = V_\nu^* \sigma$$

We obtain $\sigma = \frac{\nu-3}{2}$ or $\sigma = -1$ (which is formally included for $\nu = 1$) and $\delta = -1 + \frac{\mu\nu+\mu}{2\nu}$. The equation $V_\nu^* \delta = V_\lambda^* \sigma$ gives $\lambda(\nu+1)(\nu+\mu+\mu\nu) = \mu\nu(\nu+1)$. Since $\nu + 1 > 0$ we obtain $\lambda(\nu + \mu + \mu\nu) = \mu\nu$. The equation $V_\mu^* \delta = V_\lambda^* \delta$ leads to the same condition. An example is given by $\lambda = \frac{1}{4}, \mu = \frac{1}{2}, \nu = 1$.

## 3. Closing remarks

The question remains what can be said about a suitable dual system $(B^*, T^*)$ (see [4] for the general notion of a dual system) if no natural dual system exists. The following **conjecture** seems reasonable. There is a closed set $B^*$ such that $(B^*, T^*)$ is a fibred system where $T^*$ is piecewise defined by the transposed matrices. Generally, the set $B^*$ supports a 1-conformal measure (see e.g. [2] and [6]). This property is based on the fact that the equation

$$\sum_{k_1,\ldots,k_s} \frac{\omega^*(k_1,\ldots,k_s;y)}{1 + V^*(k_1,\ldots,k_s)y} = \frac{1}{1+y}$$

holds. In the case of a dual Moebius system $B^*$ is an interval and this measure clearly is Lebesgue measure. However, there exist examples where $B^*$ is a union of infinitely many intervals. We give one such example.

$B = [0, 1]$ and

$$Tx = \begin{cases} \frac{x}{1-2x}, & \text{if } 0 \leq x < \frac{1}{3}, \\ \frac{1-2x}{x}, & \text{if } \frac{1}{3} \leq x < \frac{1}{2}, \\ \frac{1-x}{x}, & \text{if } \frac{1}{2} \leq x < 1, \end{cases}$$

$B^* = \bigcup_{k=0}^\infty ]2k, 2k+1]$ and

$$Ty = \begin{cases} -y + 2, & \text{if } y \in \bigcup_{k=1}^\infty ]2k, 2k+1], \\ \frac{1-y}{y}, & \text{if } y \in \bigcup_{k=0}^\infty ]\frac{1}{2k+2}, \frac{1}{2k+1}], \\ \frac{1-2y}{y}, & \text{if } y \in \bigcup_{k=0}^\infty ]\frac{1}{2k+3}, \frac{1}{2k+2}]. \end{cases}$$

Therefore the invariant density for $T$ is given as

$$h(x) = \sum_{k=0}^{\infty} \left(\frac{2k+1}{1+(2k+1)x} - \frac{2k}{1+2kx}\right).$$

I want to express my sincere thanks to the referee whose advice was very helpful.